\documentclass[12pt,a4paper]{article}

\usepackage{shuffle}
\usepackage{authblk}

\usepackage{amsmath}
\usepackage{amssymb}
\usepackage{graphicx}
\usepackage{hyperref}

\usepackage[english]{babel}

\newcommand{\R}{{\mathbb R}}
\newcommand{\N}{{\mathbb N}}

\newcommand {\LL}{{\mathcal L}}
\newcommand {\PP}{{\mathcal P}}

\newcommand{\A}{{\mathcal A}}
\newcommand{\B}{{\mathcal B}}
\newcommand{\J}{{\mathcal J}}

\newcommand{\Lin}{\mathrm{Lin}}
\newcommand{\ord}{\mathrm{ord}}
\newcommand{\sign}{\mathrm{sign}}
\newcommand{\shf}{\shuffle}
\newcommand{\ad}{\mathrm{ad}}

\newtheorem{theorem}{\rm\bf Theorem}[section]

\newtheorem{remark}[theorem]{\rm\bf Remark}
\newtheorem{definition}[theorem]{\rm\bf Definition}

\title{Implementation of the algorithm for constructing homogeneous approximations of nonlinear control systems}
\author[a,b]{Grigory SKLYAR}
\author[c,d]{Pavel BARKHAYEV}
\author[e]{Svetlana IGNATOVICH}
\author[e]{Viktor RUSAKOV}
\affil[a]{Faculty of Computer Science and Information Technology\\ 
	
	West Pomeranian University of Technology in Szczecin,
		\.{Z}o\l{}nierska 49, 71-210 Szczecin, Poland\\ 
		
		e-mail: \url{gsklyar@zut.edu.pl}
	
\:}

\affil[b]{Institute of Mathematics\\ 
	
	University of Szczecin,
	Wielkopolska 15, 70-451 Szczecin, Poland\\ 
	
	e-mail: \url{grigorij.sklyar@usz.edu.pl}

\:}

\affil[c]{Department of Mathematical Sciences\\ 
	
	Norwegian University of Science and Technology, H\o{}gskoleringen 1, 7491  Trondheim, Norway\\ 
	
	e-mail: \url{pavlo.barkhayev@gmail.com}, \url{pavloba@ntnu.no}

\:}

\affil[d]{B.Verkin Institute for Low Temperature Physics and Engineering\\
	
	Academy of Sciences of Ukraine, Nauki Ave. 47, 61103 Kharkiv, Ukraine

\:}

\affil[e]{Department of Applied Mathematics\\ 
	
	V. N. Karazin Kharkiv National University, Svobody Sq. 4, 61022 Kharkiv, Ukraine\\ 
	
	e-mail: \url{ignatovich@ukr.net}, \url{s.ignatovich@karazin.ua}

\:}

\date{}                     
\begin{document}

\maketitle

\begin{abstract}
We present a ``calculator'' for constructing a homogeneous approximation of nonlinear control systems, which is based on the algebraic approach developed by the authors in their previous papers. This approach mainly uses linear algebraic and combinatorial tools, so, it is perfectly adapted to computer realization. We describe the algorithm and discuss its capabilities and limitations. We present its implementation as a web application and show by example how this app works.

\end{abstract}

\section{Introduction}

When studying properties of dynamical systems or problems of control design, the first step commonly used in mechanics and engineering is to linearize the systems under consideration. Namely, the nonlinear system is substituted by a linear one, which is nevertheless close enough to the original one, so that the results obtained for the linear system are applicable to the original problem. A brilliant example of such a method is the Lyapunov stability theory.

However, in many cases it is impossible to choose an appropriate linear system. As a simple example, let us consider a two-input driftless control system 
\begin{equation}\label{eq:G}
\dot x=X_1(x)u_1+X_2(x)u_2,  
\end{equation}
where $X_1(x)$ and $X_2(x)$ are smooth vector fields defined in a domain of $\R^n$ with $n\ge3$. Suppose that this system is locally controllable at some point $x^0\in\R^n$, i.e., any point from a neighborhood $U(x^0)$ can be reached from $x(0)=x^0$ by a piecewise continuous control. However, its linearization, which  should be of the form
\begin{displaymath}
\dot x=b_1u_1+b_2u_2,  
\end{displaymath}
where $b_1,b_2\in \R^n$ are constant vectors, necessarily is uncontrollable. Thus, any linear system looses the most important property of the original nonlinear system, therefore, hardly can be considered as its approximation.

In such a case, one can try to use an appropriate \emph{nonlinear} system as an approximation. It is natural to expect that an approximating system, being nonlinear, has a simpler structure than the original one. Deep and long-term studies in this direction were carried out by several researches, who analyzed different problems. But despite the diversity of approaches, in the end all these authors came to the same concept of a so-called \emph{homogeneous approximation}.  We mention several important original publications  \cite{Crouch:84,Bressan:85,Stefani:85,Hermes:86,Agrachev_Gamkrelidze_Sarychev:89,Bianchini_Stefani:90,Sussmann:92}; this list is undoubtedly far from being complete. The results obtained were summarized in  \cite{Hermes:91} and \cite{Bellaiche:96}. In the above works, the differential geometry language and tools were used.

Another way to introduce a nonlinear approximation was related to the interpretation of nonlinear control systems as series of iterated integrals. This idea was originated by M. Fliess \cite{Fliess:78,Fliess:81}, who proposed to apply Chen series 
\cite{Chen} to control problems. The approach uses the tools of free algebras \cite{Reutenauer} to study control systems  \cite{Fliess:80,Kawski_Sussmann:97,Kawski:97,Isidori}. 

In the mid 90s of the last century, two of the authors of the present paper, Grigory Sklyar and Svetlana Ignatovich, considered another approximation problem, namely, approximation in the sense of time optimality \cite{Sklyar_Ignatovich:96ZAMM,Sklyar_Ignatovich:00}.  
They proposed a way for analysis of the time optimal problem for nonlinear systems based on algebraic methods \cite{Sklyar_Ignatovich:03}. The main goal was to generalize the moment approach, which  proved to be an effective tool in the linear setting \cite{Korobov_Sklyar:87,Korobov_Sklyar:91}. As a result, they gave an explicit construction of approximating systems in the sense of time optimality using the free algebras language. It turned out that in essence, from the algebraic point of view, the approximation in the sense of time optimality and the homogeneous approximation are similar concepts though related to different free algebras \cite{Sklyar_Ignatovich:02,Sklyar_Ignatovich:08,Sklyar_Ignatovich:14}. In general, algebraic approach proved to be very natural and allowed a deeper understanding of nonlinear systems
\cite{Sklyar_Ignatovich_Barkhayev:05,Sklyar_Ignatovich:07,Ignatovich:09,Ignatovich:11,Sklyar_Ignatovich:14,Sklyar_Ignatovich:21}. 

It turns out that the proposed construction, in addition to its theoretical significance, has the following advantage, which is essential for computer implementation: it mainly consists of purely linear algebraic and combinatorial steps and, therefore, promises to become a universal tool adapted to arbitrary control-affine systems. However, in all their previous papers, Grigory Sklyar and Svetlana Ignatovich did not bring the approach to practical calculations, although they perfectly understood the possibilities provided by the algorithmic nature of the method. Another author of the present paper, Pavel Barkhayev, being deeply involved in the matter, initialized this work and created the first version of the algorithm realization. When the last author, Viktor Rusakov, joined, the computer program was improved and, finally, a ``calculator'' of homogeneous approximations was obtained.  

The main goal of the paper is to present and discuss the algorithm for constructing a homogeneous approximation based on the mentioned algebraic approach, its implementation, capabilities and limitations. 

Concluding the introduction section, we illustrate possible advantages of a homogeneous approximation by a very simple example. Let us consider a stabilization problem for a rotating rigid body 
\begin{equation}\label{eq:sec1_1}
\dot x_1=\alpha_1x_2x_3+u_1, \  \dot x_2=\alpha_2x_1x_3+u_2, \  \dot x_3=\alpha_3x_1x_2,
\end{equation}
where $\alpha_3\ne0$. In the paper \cite{Zuyev:19}, time-varying feedback controls were proposed, which stabilize the equilibrium $x=0$ and can be easily found explicitly. Using a homogeneous approximation, one can simplify the proposed controls. 

\begin{figure}
\centering
\includegraphics[scale=0.85]{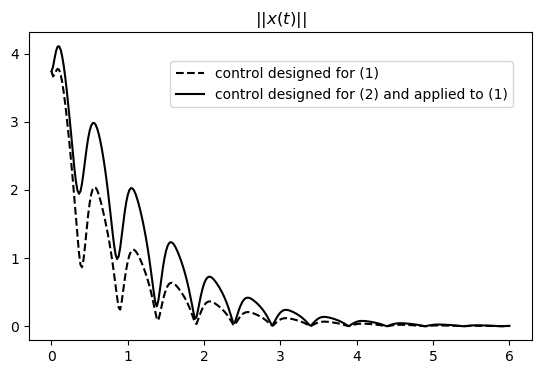}
\caption{Norms of the trajectories $x(t)$ of the system (\ref{eq:sec1_1}).}\label{Fig:f1}
\end{figure}

The system \eqref{eq:sec1_1} is not homogeneous, however, its homogeneous approximation can be easily guessed;  it has the form 
\begin{equation}\label{eq:sec1_2}
\dot x_1=u_1, \ \ \dot x_2=u_2, \ \ \dot x_3=\alpha_3x_1x_2.
\end{equation}
Now, following \cite{Zuyev:19}, we find time-varying feedback controls stabilizing the homogeneous system \eqref{eq:sec1_2}
\begin{equation}\label{eq:sec1_3}
\begin{array}{l}
u_1(t,x)=-\gamma x_1 + \frac{4\pi\sqrt{|a_{120}|}}{\varepsilon}\cos\left(\frac{2\pi t}{\varepsilon}\right), \\[3pt] 
u_2(t,x)=-\gamma x_2 + \frac{4\pi\sqrt{|a_{120}|}}{\varepsilon}\sign(a_{120})\cos\left(\frac{2\pi t}{\varepsilon}\right),  
\end{array}
\end{equation}
where $\varepsilon>0$ is rather small, $\gamma>0$, and
\begin{displaymath}
\textstyle
a_{120}= - \frac{\gamma}{2\alpha_3}x_3 - \frac12x_1x_2.
\end{displaymath}
It turns out that these controls also stabilize the original system \eqref{eq:sec1_1}. We note that they are slightly simpler than those proposed in the mentioned paper for the system \eqref{eq:sec1_1}, since they are designed for a simpler system \eqref{eq:sec1_2}. In Fig.~\ref{Fig:f1} we compare two trajectories of the system \eqref{eq:sec1_1}: the first trajectory (dotted line) corresponds to the controls proposed in \cite{Zuyev:19} and the second one (solid line) corresponds to the controls \eqref{eq:sec1_3}; in both cases  $\alpha_1=3$, $\alpha_2=2$, $\alpha_3=1$, $\gamma=5$, $\varepsilon=1$, and $x(0)=(3,2,1)^\top$. 

We note that the controls \eqref{eq:sec1_3} do not depend on $\alpha_1,\alpha_2$ and, therefore, can be applied to systems with different $\alpha_1,\alpha_2$ or with $\alpha_1,\alpha_2$ that vary within some bounds, which leads to robust control design. Moreover, the controls \eqref{eq:sec1_3} stabilize, along with \eqref{eq:sec1_1}, other systems for which \eqref{eq:sec1_2} is a homogeneous approximation.

\section{Algebraic approach to the homogeneous approximation problem}\label{sec2}

By definition, a homogeneous approximation is of a local character, so, we consider control systems in a neighborhood of a given point; without loss of generality we assume this point to be the origin. The algorithm we present in this paper is applicable in two situations: (i)~the system is autonomous, control-linear and the starting point is fixed,
\begin{equation}\label{eq:sec2_1}
\dot x = \sum_{i=1}^mX_i(x)u_i, \ x(0)=0;
\end{equation}
(ii) the system is non-autonomous, control-affine and the end point, which is an equilibrium of the system, is fixed,
\begin{equation}\label{eq:sec2_2}
\dot x=a(t,x)+\sum_{i=1}^mb_i(t,x)u_i, \ a(t,0)\equiv 0, \ x(\theta)=0.
\end{equation}

In this paper we concentrate on the second case, where we assume $m=1$ for simplicity, though, the main ideas for both cases are similar. A description of the algorithm for systems \eqref{eq:sec2_1} can be found in \cite{Sklyar_Ignatovich:20_conf}.

In the following subsections, we briefly outline our approach, which was first described in \cite{Sklyar_Ignatovich:03}.

\subsection{Series representation of control systems}

We consider a non-autonomous single-input system of the form 
\begin{equation}\label{eq:sec2_system}
\dot{x} = a(t,x)+b(t,x)u, \ \ a(t,0)\equiv 0,
\end{equation}
where $a(t,x)$, $b(t,x)$ are analytic in a neighborhood of the origin in $\R^{n+1}$, and we are interested in trajectories ending at the origin, i.e., $x(\theta)=0$, for some time moment $\theta>0$. The condition $a(t,0)\equiv0$ means that the origin is an equilibrium of the system. The set of admissible controls is defined as 
\begin{equation}\label{eq:sec2_control}
u\in B^\theta=\{u(t)\in L_\infty[0,\theta]:|u(t)|\le 1\mbox{ a.e.}\}.
\end{equation} 
Then there exists $T>0$ such that for any $\theta \in(0,T)$ and any $u\in B^\theta$ the solution of the differential equation $\dot x=a(t,x)+b(t,x)u(t)$ satisfying the end condition ${x(\theta)=0}$ is well defined for $t\in[0,\theta]$. Therefore, the system \eqref{eq:sec2_system} defines the operator $S_{a,b}$, which maps a pair $(\theta,u)$ to the starting point,
\begin{displaymath}
S_{a,b}(\theta,u)=x^0, \ \ \theta\in(0,T), \ u\in B^\theta.
\end{displaymath} 

The operator $S_{a,b}$ admits the following representation in the form of a series 
\begin{equation}\label{eq:sec2_series}
S_{a,b}(\theta, u) = \sum_{k\ge 1, m_i\ge 0} v_{m_1\ldots m_k} \xi_{m_1\ldots m_k}(\theta, u),
\end{equation}
where
\begin{equation}\label{eq:sec2_xi}
\begin{array}{l}
\xi_{m_1\ldots m_k}(\theta, u) =\\
\displaystyle = \int_0^\theta\int_0^{\tau_1}\cdots\int_0^{\tau_{k-1}}\prod_{i=1}^k(\tau_i^{m_i}u(\tau_i))d\tau_k\cdots d\tau_1
\end{array}
\end{equation}
are \emph{nonlinear power moments} and $v_{m_1\ldots m_k}\in\R^n$ are constant vector coefficients. A precise formula for $v_{m_1\ldots m_k}$ is given below, see \eqref{f_for_v}.

Since we are interested in local analysis, we assume that $\theta$ is small. Admissible controls are from the set \eqref{eq:sec2_control}, hence, the following natural asymptotic equivalence can be adopted for nonlinear power moments \eqref{eq:sec2_xi}:
\begin{equation}\label{eq:xi_equiv}
\xi_{m_1\ldots m_k}(\theta,u)\sim \theta^{m_1+\cdots+m_k+k},
\end{equation}
i.e., the number $m=m_1+\cdots+m_k+k$ can be considered as an \emph{order} of the integral $\xi_{m_1\ldots m_k}(\theta,u)$. 

This order naturally defines the homogeneity property, which allows us to introduce a homogeneous approximation of the system \eqref{eq:sec2_system} in terms of its series. 

\begin{definition}{}\label{def:ha1}
Suppose that, after some coordinate transformation in the system \eqref{eq:sec2_system}, the series $S_{a,b}(\theta,u)$ takes the form
\begin{equation}\label{eq:sys_after_tr}
(S_{a,b}(\theta,u))_k=z_k(\theta,u)+\rho_k(\theta,u), \ k=1,\ldots,n,
\end{equation}
where
\begin{itemize}
\item[(i)] $z_k(\theta,u)$ is a sum of integrals \eqref{eq:sec2_xi} of order $w_k$,
\item[(ii)] $\rho_k(\theta,u)$ contains integrals \eqref{eq:sec2_xi} of order greater than $w_k$,
\item[(iii)] the system
\begin{equation}\label{eq:sys_appr}
\dot x=\widehat a(t,x)+\widehat b(t,x)u,
\end{equation}
which corresponds to the series 
\begin{displaymath}
(S_{\widehat a,\widehat b}(\theta,u))_k=z_k(\theta,u), \ k=1,\ldots,n,
\end{displaymath}
is controllable.
\end{itemize}
Then \eqref{eq:sys_appr} is called a homogeneous approximation of the system \eqref{eq:sec2_system}.
\end{definition}

It might seem that the problem is to find a coordinate transformation taking the series to the form \eqref{eq:sys_after_tr}. However, an approximation itself is defined by $z_k(\theta,u)$, which are linear combinations of integrals \eqref{eq:sec2_xi}. The elements $z_k(\theta,u)$ can be found directly, without finding an appropriate transformation. Moreover, only basic linear algebraic and simple combinatorial tools can be used, as is shown in the next subsections.   

\subsection{Algebra of nonlinear power moments}

It turns out that the linear span of nonlinear power moments \eqref{eq:sec2_xi} can be considered as a free associative algebra 
\begin{displaymath}
\A=\Lin\{\xi_{m_1\ldots m_k}:m_1,\ldots,m_k\ge0, \ k\ge1\},
\end{displaymath}
equipped with the (formal) multiplication defined on basis elements as
\begin{displaymath}
\xi_{m_1\ldots m_k}\xi_{n_1\ldots n_r}=\xi_{m_1\ldots m_kn_1\ldots n_r}.
\end{displaymath}
The series $S_{a,b}$ defines the linear map $v:\A\to\R^n$ by
\begin{displaymath}
v(\xi_{m_1\ldots m_k})=v_{m_1\ldots m_k}.
\end{displaymath}
The asymptotic property \eqref{eq:xi_equiv} leads to a \emph{canonical grading structure} in the algebra $\A$ defined as 
\begin{equation}\label{eq:sec2_grading_canonical}
{\mathcal A} = \sum_{m=1}^\infty \A^m,
\end{equation}
where
\begin{displaymath}
\A^m=\Lin\{\xi_{m_1\ldots m_k}:m_1+\cdots+m_k+k=m\}
\end{displaymath}
are the subspaces of homogeneous elements of order $m$. Below we write $\ord(z)=m$ if $z\in\A^m$. We emphasize that subspaces $\A^m$ are finite-dimensional. 

The free associative algebra $\A$ is closely related to the free Lie algebra generated by the elements $\xi_m$,  $m\ge0$, with the Lie bracket operation $[\ell_1,\ell_2]=\ell_1\ell_2-\ell_2\ell_1$; we denote this free Lie algebra by $\LL$. The Lie algebra $\LL$ is graded as well, 
\begin{displaymath}
\LL=\sum_{m=1}^\infty\LL^m, \ \ \ \LL^m=\LL\cap\A^m.
\end{displaymath}
Omitting details, which can be found in
\cite{Sklyar_Ignatovich:14},
we mention two properties of the map~$v$. 

1) Below we assume that the system \eqref{eq:sec2_system} is locally accessible, i.e., the Rashevsky-Chow condition is satisfied. This can be expressed as
\begin{equation}\label{eq:sec2_Rash}
v(\LL)=\R^n.
\end{equation}

2) The series \eqref{eq:sec2_series} is ``realizable'', i.e., it corresponds to a control system. In our case this requirement reduces to the following property
\begin{equation}\label{eq:sec2_realiz}
\mbox{if }  \ell\in\LL \mbox{ and } v(\ell)=0, \mbox{ then }  v(\ell z)=0 \mbox{ for any } z\in\A.
\end{equation}

An analytic change of coordinates in the system \eqref{eq:sec2_system} can be expressed as a transformation over the series $S_{a,b}$. In order to find a result of such a transformation, one needs to find products of integrals \eqref{eq:sec2_xi}. The product of two integrals corresponds to a \emph{shuffle} operation in the algebra $\A$, which can be defined recursively as 
\begin{displaymath}
	\xi_{m_1}\shf\xi_{n_1}=\xi_{m_1n_1}+\xi_{n_1m_1},
\end{displaymath}
\begin{displaymath}
	\xi_{m_1}\shf\xi_{n_1\ldots n_r}=\xi_{m_1n_1\ldots n_r}+\xi_{n_1}(\xi_{m_1}\shf\xi_{n_2\ldots n_r}),
\end{displaymath}
\begin{displaymath}
\xi_{m_1\ldots m_k}\shf\xi_{n_1\ldots n_r}=
\end{displaymath}
\begin{displaymath}
=\xi_{m_1}(\xi_{m_2\ldots m_k}\shf\xi_{n_1\ldots n_r})+\xi_{n_1}(\xi_{m_1\ldots m_k}\shf\xi_{n_2\ldots n_r}).
\end{displaymath}
For example, the product
\begin{displaymath}
\xi_{m_1}(u,\theta)\cdot \xi_{n_1}(u,\theta)=\int_0^\theta\tau_1^{m_1}u(\tau_1)d\tau_1\cdot\int_0^\theta\tau_2^{n_1}u(\tau_2)d\tau_2
\end{displaymath}
equals the integral of the function $\tau_1^{m_1}\tau_2^{n_1}u(\tau_1)u(\tau_2)$ over the square $\{(\tau_1,\tau_2):0\le\tau_1,\tau_2\le\theta\}$, which can be written out as a sum of two integrals of the same function over the triangles $\{(\tau_1,\tau_2):0\le\tau_1\le\tau_2\le\theta\}$ and $\{(\tau_1,\tau_2):0\le\tau_2\le\tau_1\le\theta\}$. Therefore,
\begin{displaymath}
\xi_{m_1}(\theta,u)\cdot \xi_{n_1}(\theta,u)=\xi_{m_1n_1}(\theta,u)+\xi_{n_1m_1}(\theta,u),
\end{displaymath}
which corresponds to the shuffle product in $\A$.

\subsection{Finding a series representation for a homogeneous approximation} 

Now let us return to Definition~\ref{def:ha1}. Justifying the terms ``homogeneous'' and ``approximation'', it also gives a clue to finding a homogeneous approximation. Actually, one can start with the series $S_{a,b}$ and transform it, aiming to separate the ``main part'' with respect to the order in $\A$. The goal is to find a set of homogeneous elements $z_1,\ldots,z_n$ such that the ``series'' $S=(z_1,\ldots,z_n)^\top$ is realizable and accessible, i.e., satisfied conditions 1) and 2) mentioned above (see \eqref{eq:sec2_Rash} and \eqref{eq:sec2_realiz}). Desired elements can be obtained step by step; many examples can be found in 
\cite{Sklyar_Ignatovich_Barkhayev:05}.

However, one can obtain elements $z_1,\ldots,z_n$ explicitly, without such a brutal force way. For a given series $S_{a,b}$ of the form \eqref{eq:sec2_series} or, what is the same, for the linear map $v$ satisfying conditions 1) and 2), we introduce the linear subspaces
\begin{displaymath}
\PP^m=\{\ell\in\LL^m:v(\ell)\in v(\LL^1+\cdots+\LL^{m-1})\}, \ m\ge1,
\end{displaymath}
and denote
\begin{displaymath}
\LL_{a,b}=\sum_{m=1}^\infty\PP^m.
\end{displaymath}
It can be shown that $\LL_{a,b}$ is a Lie subalgebra of $\LL$; by definition, it is graded and has codimension $n$ in $\LL$. We call it \emph{a core Lie subalgebra} of the system \eqref{eq:sec2_system}. Also we introduce the following graded \emph{right ideal} generated by $\LL_{a,b}$,
\begin{displaymath}
\J_{a,b}=\Lin\{\ell z:\ell\in\LL_{a,b}, z\in\A+\R\}.
\end{displaymath}
Finally, we introduce the inner product $\langle\cdot,\cdot\rangle$ in the algebra $\A$ declaring elements $\xi_{m_1\ldots m_k}$ to form an orthonormal basis. Below we denote by $\widetilde{z}$ the orthogonal projection of an element $z$ onto the subspace $\J_{a,b}^\perp$, i.e., onto the orthogonal complement to the subspace  $\J_{a,b}$. 

\begin{theorem}{}
For a given system \eqref{eq:sec2_system}, let us choose $n$ homogeneous elements $\ell_1,\ldots,\ell_n\in\LL$ such that
\begin{displaymath}
\Lin\{\ell_1,\ldots,\ell_n\}+\LL_{a,b}=\LL.
\end{displaymath}
Then the elements $z_1,\ldots,z_n$ from Definition~\ref{def:ha1} can be chosen as $z_k=\widetilde{\ell}_k$, $k=1,\ldots,n$. Moreover, all possible sets of such elements can be obtained from the mentioned elements by polynomial transformations that preserve the order. 
\end{theorem}

For the proof, examples and comments, see
\cite{Sklyar_Ignatovich:03}, \cite{Sklyar_Ignatovich_Barkhayev:05}.

Thus, the homogeneous approximation is, in essence, unique and is defined by the core Lie subalgebra $\LL_{a,b}$. Moreover, finding a series $S_{\widehat a,\widehat b}$ for the homogeneous approximation can be carried out in a purely algebraic way. 

\subsection{Constructing an approximating system: general case}\label{sec24}

First, we notice that the system \eqref{eq:sys_appr} is not unique. However, if $\widehat a(t,x)\equiv0$, then $\widehat b(t,x)$ is defined uniquely. We can construct components $\widehat b_i(t,x)$ of the vector function $\widehat b(t,x)=(\widehat b_1(t,x),\ldots,\widehat b_n(t,x))^\top$ one by one in the following way.

According to Definition~\ref{def:ha1}, the element $z_k=\widetilde{\ell}_k$ is of order $w_k$, $k=1,\ldots,n$. To start with, one can show that $\widetilde\ell_1=\alpha_1\xi_{w_1-1}$, where $\alpha_1\in\R$ is nonzero. Then we set $\widehat{b}_1(t,x)=-\alpha_1 t^{w_1-1}$. Further we act by induction. For any $i=2,\ldots,n$, we write the element $z_i=\widetilde{\ell}_i$ in the form
\begin{displaymath}
\widetilde\ell_i=\sum_{j=0}^{w_i-2}y_j\xi_j+\alpha_i\xi_{w_i-1}, \ y_j\in\A^{w_i-j-1}, \ \alpha_i\in\R.
\end{displaymath}
One can show that $y_j$ are shuffle polynomials of the elements $\widetilde\ell_1,\ldots,\widetilde{\ell}_{i-1}$. This follows form the fact that $y_j\in\J_{a,b}^\perp$ and elements $\widetilde\ell_1^{\shf q_1}\shf\cdots\shf\widetilde\ell_n^{\shf q_n}$ form a basis of the subspace $\J_{a,b}^\perp$. Here and below we use the notation for the ``shuffle power'': $\ell^{\shf q}=\ell\shf\cdots\shf\ell$ ($q$ times). Therefore,
\begin{displaymath}
y_j=p_j(\widetilde\ell_1,\ldots,\widetilde{\ell}_{i-1}), \ 0\le j\le w_i-2,
\end{displaymath}
where $p_j$ denotes a shuffle polynomial,
\begin{displaymath}
p_j(\widetilde\ell_1,\ldots,\widetilde{\ell}_{i-1})=\sum\beta^j_{q_1\ldots q_{i-1}}\widetilde\ell_1^{\shf q_1}\shf\cdots\shf\widetilde\ell_{i-1}^{\shf q_{i-1}}
\end{displaymath}
and the sum is taken over all integers $q_1,\ldots,q_{i-1}\ge0$ such that $w_1q_1+\cdots+w_{i-1}q_{i-1}=w_i-j-1$. Then we set
\begin{displaymath}
\widehat{b}_i(t,x)=-\sum_{j=0}^{w_i-2}P_j(x_1,\ldots,x_{i-1})t^j-\alpha_i t^{w_i-1},
\end{displaymath}
where
\begin{displaymath}
P_j(x_1,\ldots,x_{i-1})=\sum\beta^j_{q_1\ldots q_{i-1}}x_1^{q_1}\cdots x_{i-1}^{q_{i-1}}.
\end{displaymath}
So, the polynomial $P_j$ in real variables $x_1,\ldots,x_{i-1}$ is constructed as follows: for any monomial $x_1^{q_1}\cdots x_{i-1}^{q_{i-1}}$, its coefficient is the same as the coefficient of the monomial $\widetilde\ell_1^{\shf q_1}\shf\cdots\shf\widetilde\ell_{i-1}^{\shf q_{i-1}}$ in the shuffle polynomial $p_j$.

\subsection{Constructing an approximating system: autonomous case}\label{sec25}

If the initial system is autonomous, i.e., has the form 
\begin{equation}\label{eq:sec2_syst_aut}
\dot x=a(x)+b(x)u, \ \ a(0)=0,
\end{equation}
then a homogeneous approximation can also be chosen as an autonomous system, which is uniquely defined. It can be constructed as follows.

Let us introduce a \emph{derivation} $\varphi$ in the algebra $\A$ defined by 
\begin{displaymath}
\varphi(\xi_0)=0, \  \varphi(\xi_{m})=m\xi_{m-1}, \ m\ge1.
\end{displaymath}
Recall that a derivation in the algebra is a linear map satisfying the Leibniz product rule. Hence, 
\begin{displaymath}
\varphi(\xi_{m_1\ldots m_k})=\sum_{i:m_i\ge1}m_i\xi_{m_1\ldots(m_i-1)\ldots m_k}.
\end{displaymath}
One can show that for the autonomous system \eqref{eq:sec2_syst_aut}, if $z\in\J_{a,b}^\perp$, then $\varphi(z)\in\J_{a,b}^\perp$. 

Denote by $\psi$ a linear map in $\A$ defined as 
\begin{displaymath}
\psi(\xi_0)=1, \ \ \psi(\xi_{m_1\ldots m_{k-1}0})=\xi_{m_1\ldots m_{k-1}}
\end{displaymath}
and $\psi(\xi_{m_1\ldots m_k})=0$ if $m_k\ge1$. Again, if $z\in\J_{a,b}^\perp$ and $\ord(z)\ge2$, then $\psi(z)\in\J_{a,b}^\perp$. 

We notice that if $z\in\A^m$, then $\varphi(z)\in\A^{m-1}$ and $\psi(z)\in\A^{m-1}$ if $m\ge2$.

Now we are ready to construct the autonomous homogeneous approximation 
\begin{displaymath}
\dot x=\widehat a(x)+\widehat b(x)u, \ \ \widehat a(0)=0,
\end{displaymath}
for the system \eqref{eq:sec2_syst_aut}. First, one can show that $\widetilde\ell_1=\alpha_1\xi_{0}$, where $\alpha_1\ne0$; we set $\widehat{a}_1(x)=0$, $\widehat{b}_1(x)=-\alpha_1$. Then we act by induction. For any $i=2,\ldots,n$, we consider the elements $\varphi(\widetilde\ell_i)$ and $\psi(\widetilde\ell_i)$; since they belong to $\J_{a,b}^\perp$, they can be written in the form 
\begin{displaymath}
\varphi(\widetilde\ell_i)=p_{1,i}(\widetilde\ell_1,\ldots,\widetilde\ell_{i-1}), \ \psi(\widetilde\ell_i)=p_{2,i}(\widetilde\ell_1,\ldots,\widetilde\ell_{i-1}),
\end{displaymath}
where $p_{1,i}$ and $p_{2,i}$ are shuffle polynomials,
\begin{displaymath}
p_{j,i}(\widetilde\ell_1,\ldots,\widetilde{\ell}_{i-1})=\sum\beta^j_{q_1\ldots q_{i-1}}\widetilde\ell_1^{\shf q_1}\shf\cdots\shf\widetilde\ell_{i-1}^{\shf q_{i-1}},
\end{displaymath}
$j=1,2$, and the sums are taken over all integers $q_1,\ldots,q_{i-1}\ge0$ such that $w_1q_1+\cdots+w_{i-1}q_{i-1}=w_i-1$. Then we set
\begin{displaymath}
\begin{array}{c}
\widehat{a}_i(x)=-P_{1,i}(x_1,\ldots,x_{i-1}), \\ 
\widehat{b}_i(x)=-P_{2,i}(x_1,\ldots,x_{i-1}),
\end{array}
\end{displaymath}
where
\begin{displaymath}
P_{j,i}(x_1,\ldots,x_{i-1})=\sum\beta^j_{q_1\ldots q_{i-1}}x_1^{q_1}\cdots x_{i-1}^{q_{i-1}}.
\end{displaymath}
Thus, the polynomial $P_{j,i}$ in real variables $x_1,\ldots,x_{i-1}$ is constructed so that for any monomial $x_1^{q_1}\cdots x_{i-1}^{q_{i-1}}$, its coefficient is the same as the coefficient of the monomial $\widetilde\ell_1^{\shf q_1}\shf\cdots\shf\widetilde\ell_{i-1}^{\shf q_{i-1}}$ in the shuffle polynomial $p_{j,i}$.

\section{The algorithm}\label{sec3}

The pipeline of the algorithm for constructing a homogeneous approximation for a given affine control system may be described naturally in a step-by-step manner. In what follows we give a detailed description of the algorithm splitted into nine steps.

We illustrate all the steps by the following example:
\begin{equation}\label{eq_ex1}
\left\{
\begin{array}{l}
\dot{x}_1 = - u \cos x_1,\\
\dot{x}_2 = -\sin^2 x_1 + t^2 u, \\
\dot{x}_3 = 2x_1^2 \sin t - x_2 u,
\end{array}
\right.
\end{equation}
i.e., for $n=3$, $a(t,x) = (0, -\sin^2 x_1, 2x_1^2 \sin t)^T$, $b(t,x) = (-\cos x_1, t^2, -x_2)^T$.

\smallskip

\noindent\textbf{Step 1.} Construction of the algebra $\A$.

We fix $N\in\N$ and for any $m=1,\ldots,N$ we construct the graded subspace $\A^m$, i.e., we find all the sequences $\{m_1,\ldots,m_k\}$ of length $k\le m$ satisfying the condition $m_1+\cdots + m_k=m-k$ (admissible sequences). The corresponding nonlinear power moments $\xi_{m_1\ldots m_k}$ belong to $\A^m$.
 
To find these sequences, one may consider all ${(m-k+1)}$-ary numbers with no more than $k$ digits and check the sum of the digits. We note that $\dim \A^m = 2^{m-1}$.
 
Having in mind the system \eqref{eq_ex1}, in Table~\ref{tab:1} we give the admissible sequences for $N=4$.

\begin{table}[!h]
	\caption{Basis of the algebra $\A$.}\label{tab:1}
	\center\begin{tabular}{|c|l|}
		\hline
		$m$ & Admissible sequences \\
		\hline\hline
		1 & \{0\}  \\
		2 & \{1\}; \{0, 0\}\\
		3 & \{2\}; \{0, 1\}; \{1, 0\}; \{0, 0, 0\} \\
		4 & \{3\}; \{0, 2\}; \{2, 0\}; \{1, 1\}; \{0, 0, 1\}; \{0, 1, 0\}; \\
		& \{1, 0, 0\}; \{0, 0, 0, 0\} \\
		\hline
	\end{tabular}
\end{table}
 
\smallskip

\noindent\textbf{Step 2.} Construction of the Lie algebra $\LL$. 

Lie elements of the form 
\begin{equation}\label{eq_dynkin_bracket}
[\xi_{m_1},[\ldots [\xi_{m_{k-1}}, \xi_{m_k}]\ldots]], \quad k\ge 1, m_i\ge 0
\end{equation}
are called right normed elements \cite{Reutenauer,Kawski:01}.
It is known that the set of all right normed elements is complete in the Lie algebra $\LL$ and does not form a basis. However, discarding linearly dependent elements from the set of right normed elements, one can construct a basis of $\LL$.

Along with the canonical grading \eqref{eq:sec2_grading_canonical} in $\A$, we consider the \emph{length grading}
\begin{displaymath}
\A = \sum_{k=1}^\infty \B^k,
\end{displaymath}
where $\B^k=\Lin\{\xi_{m_1\ldots m_k}:m_1,\ldots,m_k\ge0\}$. We note that every right normed element is homogeneous in the sense that it belongs to some subspace $\A^m$ as well as to some subspace $\B^k$. This means that a basis of $\LL$ may be constructed independently in every single subspace ${\A^m\cap\B^k}$, $m, k \ge 1$. 

For any $m=1,\ldots,N$ (the number $N$ has been chosen at Step 1) and $k\le m$ we consider admissible sequences $\{m_1,\ldots,m_k\}$ such that $\xi_{m_1\ldots m_k}\in \A^m\cap\B^k$ and construct the corresponding right normed elements \eqref{eq_dynkin_bracket}.

To find linearly independent elements, we consider a realization of the subspaces $\A^m\cap\B^k$ in $\R^q$, $q=\dim(\A^m\cap\B^k)$. Namely, let $F_{m,k}$ be a bijection from the set of all elements $\xi_{m_1\ldots m_k}\in \A^m\cap\B^k$ to the set of unit vectors $\{e_j\}_{j=1}^q\in \R^q$, and we extend it by linearity to the mapping from $\A^m\cap\B^k$ to $\R^q$. Further, we represent right normed elements from $\A^m\cap\B^k$ as sums of $\xi_{m_1\ldots m_k}$ and apply $F_{m,k}$.
The obtained vectors in $\R^q$ are being checked for linear independence (adding one by one), and we discard linearly dependent vectors. 

We store the basis of the Lie algebra $\LL$ as the set of sequences $\{m_1,\ldots,m_k\}$ and as the corresponding vectors in $\R^q$.

Further we refer to the obtained set of linearly independent right normed elements as $\{g_i\}_{i\ge 1}$, assuming that $\ord(g_i) \le \ord(g_j)$ as $i<j$.

Table~\ref{tab:2} gives elements $g_i$ for $N=4$.

\begin{table}[!h]
	\caption{Basis of the Lie algebra $\LL$.}\label{tab:2}
	\center\begin{tabular}{|c|l|}
		\hline
		$m$ & $g_i$ \\
		\hline\hline
		1 & $g_1 = \xi_{0}$  \\
		2 & $g_2 = \xi_{1}$\\
		3 & $g_3 = \xi_{2}$, $g_4 = [\xi_{0}, \xi_{1}]$ \\
		4 & $g_5 = \xi_{3}$, $g_6 = [\xi_{0}, \xi_{2}]$, $g_7 = [\xi_{0}, [\xi_{1}, \xi_{0}]]$  \\
		\hline
	\end{tabular}
\end{table}

\begin{remark}{}
We note that there exist purely algebraic methods of constructing a basis of the Lie algebra: Hall, Lyndon, Shirshov bases. All these algorithms are based on the remarkable Lazard elimination theorem, see, e.g., \cite{Reutenauer}.
\end{remark}

\begin{remark}{}
Dimensions of the subspaces $\LL^m$ are be given explicitly by  Witt's formula
\begin{displaymath}
	\dim \LL^m =\frac{1}{m}\mathop\sum_{d|m}\mu(d)2^{m/d} \ \mbox{ for } m>1,
\end{displaymath}
where  $\mu: \N\rightarrow \{-1,0,1\}$ is a M\"obius function,
\begin{displaymath}
	\mu(d) = \left\{
	\begin{array}{ll}
		1, & d \ \mbox{ is a square-free positive integer with} \\
		& \mbox{an even number of prime factors;}\\
		-1, & d \ \mbox{ is a square-free positive integer with} \\
		& \mbox{an odd number of prime factors;}\\
		0, & d \ \mbox{ has a squared prime factor.}
	\end{array}
	\right.
\end{displaymath}
In Table~\ref{tab:3} we give comparison between $\dim \A^m$ and $\dim \LL^m$ for several first values of $m$.

\begin{table}[!h]
	\caption{Dimensions of $\A^m$ and $\LL^m$.}\label{tab:3}
	\setlength\tabcolsep{4pt}
	\center\small\begin{tabular}{|c||c|c|c|c|c|c|c|c|c|c|}
		\hline
		$m$ & 1 & 2 & 3 & 4 & 5 & 6 & 7 & 8 & 9 & 10 \\
		\hline
		$\dim \A^m$ & 1 & 2 & 4 & 8 & 16 & 32 & 64 & 128 & 256 & 512 \\
		\hline
		$\dim \LL^m$ & 1 & 1 & 2 & 3 & 6 & 9 & 18 & 30 & 56 & 99 \\
		\hline
	\end{tabular}
\end{table}
\end{remark}

\begin{remark}{} 
Steps 1 and 2 do not depend on a specific control system; here we construct and store basis elements of the free associative algebra $\A$ and the corresponding Lie algebra $\LL$, which can be used later for any specific system.
	
We also note that both steps allow parallelizing in computing.
\end{remark}

\smallskip

\noindent\textbf{Step 3.} Computation of vector coefficients $v_{m_1\ldots m_k}$ of the series \eqref{eq:sec2_series}.

Let us introduce differential operators $R_a$, $R_b$ acting on a function $\varphi=\varphi(t,x)$ as
\begin{displaymath}
R_a \varphi =\varphi_t + \varphi_x a,\quad  R_{b} \varphi =\varphi_x b,
\end{displaymath}
where $\varphi_t=\frac{\partial \varphi}{\partial t}$ and $\varphi_x=(\frac{\partial \varphi}{\partial x_1},\ldots,\frac{\partial \varphi}{\partial x_n})$. Below we apply $R_a$, $R_b$ to vector functions assuming that they act componentwise. Also, let us introduce the iterative operators $\ad_{R_a}^j$ acting on $R_b$ as
\begin{displaymath}
\ad_{R_a}^{0} R_{b}=R_{b}, \quad \ad_{R_a}^j R_{b}=[R_{a}, \ad_{R_a}^{j-1}R_{b}], \ j\ge 1,
\end{displaymath}
where the bracket means the operator commutator.
Then coefficients $v_{m_1\ldots m_k}$ of the series \eqref{eq:sec2_series} equal
\begin{equation}\label{f_for_v}
v_{m_1\ldots m_k}=\frac{(-1)^{k}}{m_1!\cdots m_k!}
\ad_{R_a}^{m_1}R_{b}\cdots 
\ad_{R_a}^{m_k}R_{b}E(x){\phantom{*}\atop\Bigl|_{{t=0\atop x=0}\atop\phantom{x}}},
\end{equation}
where $E(x)\equiv x$ (see, e.g., \cite{Sklyar_Ignatovich:03}).

We implement the differential operators $R_a$, $R_b$ by using symbolic differentiation by $x$ and $t$, and the operators ${\rm ad}^j_{R_a}$ are implemented recursively. We note that the usage of recursion may be avoided, however in this case we would have to store a set of functions.

For any $m=1,\ldots,N$ and any $\xi_{m_1\ldots m_k}\in \A^m$ we apply the operators $\ad^{m_i}_{R_a}R_b$ iteratively to obtain $v_{m_1\ldots m_k}$. 

For example, for the system \eqref{eq_ex1}, finding coefficients $v_{m_1\ldots m_k}$ for all admissible sequences $\{m_1,\ldots, m_k\}$ (see Table~\ref{tab:1}), we obtain the first terms of the (infinite) series~\eqref{eq:sec2_series}:
\begin{displaymath}
S_{a,b}=\begin{pmatrix}
1 \\0 \\0 
\end{pmatrix}
\xi_{0}
- 
\begin{pmatrix}
0 \\1 \\0 
\end{pmatrix}
\xi_{2}
+ 
\begin{pmatrix}
0 \\2 \\0 
\end{pmatrix}
\xi_{01}
- 
\begin{pmatrix}
1 \\0 \\0 
\end{pmatrix}
\xi_{000}
\end{displaymath}
\begin{displaymath}
- 
\begin{pmatrix}
0 \\0 \\ 1 
\end{pmatrix}
\xi_{20}
- 
\begin{pmatrix}
0 \\0 \\2 
\end{pmatrix}
\xi_{02}
+ 
\begin{pmatrix}
0 \\0 \\2 
\end{pmatrix}
\xi_{010}
- 
\begin{pmatrix}
0 \\0 \\ 2 
\end{pmatrix}
\xi_{001}
+\ldots
\end{displaymath}

\smallskip

\noindent\textbf{Step 4.} Computation of coefficients of basis Lie elements $v(g_i)$.

Here we use the representation of the basis Lie elements in $\R^q$ as linear combinations of basis elements $\xi_{m_1\ldots m_k}$ and then use the values of $v_{m_1\ldots m_k}$ computed at the previous step taking into account the linearity of the map $v$.

For our example (\ref{eq_ex1}), using Table~\ref{tab:2} and taking into account that
\begin{displaymath}
\begin{array}{l}
g_4=[\xi_0,\xi_1]=\xi_{01}-\xi_{10}, \ \
g_6=[\xi_0,\xi_2]=\xi_{02}-\xi_{20},\\
g_7=[\xi_0,[\xi_1,\xi_0]]=2\xi_{010}-\xi_{001}-\xi_{100},
\end{array}
\end{displaymath}
we have
\begin{displaymath}
v(g_{1})=\begin{pmatrix}
	1 \\0 \\0 
\end{pmatrix}, \ v(g_{2})=v(g_{5})=0,  \ v(g_{3})=\begin{pmatrix}
0 \\-1 \\0 
\end{pmatrix},  
\end{displaymath}
\begin{displaymath}
v(g_{4})=\begin{pmatrix}
0 \\2 \\0 
\end{pmatrix}, \ v(g_{6})=\begin{pmatrix}
0 \\0 \\-1 
\end{pmatrix}, \ v(g_{7})=\begin{pmatrix}
0 \\0 \\6 
\end{pmatrix}.
\end{displaymath}

\smallskip

\noindent\textbf{Step 5.} Finding the first (by order) $n$ Lie elements $\{\ell_{j}\}_{j=1}^n$ with linearly independent vector coefficients and the generating set of the right ideal $\J_{a,b}$.

Successively for  $m=1,\ldots,N$ we divide the basis Lie elements $\{g_i\}\subset \LL^m$ into two groups as follows.

At the beginning, both groups are empty and we put $g_1$ to the first group if $v(g_1)\not=0$ and to the second one otherwise. In the latter case we proceed until we find $g_i$ such that $v(g_i)\not=0$ and set it to the first group. The first group is non-empty now. 

General step: suppose there are $k$ ($k<n$) elements in the first group denoted by 
$\ell_{1}, \ldots, \ell_{k}$, and we examine the basis Lie element $g_i$ with $\ord(g_i)=m$.

If $v(g_i)\not\in\Lin\{v(\ell_{1}), \ldots, v(\ell_{k})\}$, then we add $g_i$ to the first group and denote  $\ell_{k+1}=g_i$.

If not, then, if $v(g_i)\in v(\LL^1+\cdots+\LL^{m-1})$, we add $g_i$ to the second group. 

If not again, we find a linear combination of elements $\ell_j,\ldots,\ell_k\in\LL^m$ such that $v(g_i + c_j \ell_{j} + \cdots + c_k \ell_{k})\in v(\LL^1+\cdots+\LL^{m-1})$ and add this linear combination to the second group.

We perform this splitting until the first group consists of $n$ elements and then until all $g_i$ such that $\ord(g_i) = \ord(\ell_n)$ are being checked.

We denote the elements of the second group by $\{d_j\}$; below we use them to construct the right ideal $\J_{a,b}$.

For our example \eqref{eq_ex1} we get
\begin{displaymath}
\begin{array}{lll}
\ell_{1} = g_{1}, & \ell_{2} = g_{3}, & \ell_{3} = g_{6}
\end{array}
\end{displaymath}
and 
\begin{displaymath}
\begin{array}{llll}
d_{1} = g_{2}, & d_{2} = g_{4} + 2 g_{3}, & d_{3} = g_{5}, & d_{4} = g_{7} + 6g_{6}.
\end{array}
\end{displaymath}
Let us explain the last step, i.e., the choice of $d_4$. For $g_7\in\LL^4$ we see that $v(g_7)\in\Lin\{v(\ell_1),v(\ell_2),v(\ell_3)\}$, however,  $v(g_7)\not\in v(\LL^1+\LL^2+\LL^3)=\Lin\{v(\ell_1),v(\ell_2)\}$. Then we find that $v(g_7+6\ell_{3})\in \Lin\{v(\ell_1),v(\ell_2)\}$ (more specifically, $v(g_7+6\ell_{3})=0$), therefore, we set $d_{4} =g_{7} + 6\ell_3= g_{7} + 6g_{6}$.
\smallskip

\noindent\textbf{Step 6.} Construction of the right ideal $\J_{a,b}$. 

For every $m\in\{\ord(\ell_i): i=1,\ldots,n\}$ we construct the intersection $\J_{a,b}\cap \A^m$ of the right ideal with the subspace of elements of order $m$.
To this end, we take every $d_j$ such that $\ord(d_j) \le m$ and construct the set 
$$
D_j^m = \{d_j \xi_{s_1\ldots s_k}: \; \xi_{s_1\ldots s_k}\in \A^r, \ r=m-\ord(d_j)\}
$$
if $\ord(d_j) < m$ and $D_j^m =\{d_j\}$ if $\ord(d_j) = m$.
Due to construction one has
$$
\J_{a,b}\cap \A^m = 
\Lin \{D_j^m: \ord(d_j) \le m\}.
$$
Further we find the representation of the obtained elements $\{D_j^m: \ord(d_j) \le m\}$ in $\R^q$, $q=\dim \A^m$, considering them as linear combinations of basis elements $\xi_{m_1\ldots m_k}$. Similarly to Step~1, we fix a bijection $F_m$ from the set of all elements $\xi_{m_1\ldots m_k}\in \A^m$ to the set of unit vectors $\{e_i\}_{i=1}^q\in\R^q$, $q=\dim\A^m$. Then, for any element $d_j \xi_{s_1\ldots s_k}\in D_j^m$, we express $F_m(d_j \xi_{s_1\ldots s_k})=c_1e_1+\cdots+c_qe_q$ and compose a matrix $J_m$ with rows $(c_1,\ldots,c_q)$. 

It may happen that the rank of the matrix $J_m$ is not full. 
E.g., suppose that $\xi_0,\xi_1\in \J_{a,b}$, then necessarily $[\xi_0,\xi_1] \in \J_{a,b}$. If at the previous step we denoted $d_1=\xi_0$, $d_2=\xi_1$, $d_3=[\xi_0,\xi_1]$, then $d_3\in D_3^3$ is a linear combination of elements from $D_1^3$ and $D_2^3$ (namely, of $d_1\xi_1$ and $d_2\xi_0$). In this case the row of $J_3$ corresponding to $d_3$ linearly depends on the rows corresponding to $d_1\xi_1$ and $d_2\xi_0$.

It is crucial for the next step that $J_m$ is of full rank, so we discard linearly dependent rows from it.

For the system \eqref{eq_ex1}, $\{\ord(\ell_i): i=1,\ldots,n\} = \{1,3,4\}$,  and we get
\begin{displaymath}
\begin{array}{rcl}
\J_{a,b}\cap \A^1 & = & \{0\}, \\
\J_{a,b}\cap \A^3 & = & \Lin\{D_1^3, D_2^3\},\\
\J_{a,b}\cap \A^4 & = & \Lin\{D_1^4, D_2^4, D_3^4, D_4^4\},
\end{array}
\end{displaymath}
where
\begin{displaymath}
\begin{array}{rcl}
D_1^3 & = & \{d_1 \xi_0\} = \{\xi_{10}\}, \\
D_2^3 & = & \{d_2\} = \{\xi_{01} - \xi_{10} + 2\xi_{2} \},\\
D_1^4 & = & \{d_1 \xi_{00}, d_1 \xi_1\} = \{\xi_{100}, \xi_{11}\},\\
D_2^4 & = & \{d_2 \xi_0\} = \{\xi_{010} - \xi_{100} + 2\xi_{20}\},\\
D_3^4 & = & \{d_3\} = \{\xi_{3} \},\\
D_4^4 & = & \{d_4\} = \{2\xi_{010}-\xi_{001} -   \xi_{100} + 6\xi_{02} - 6\xi_{20}\}.
\end{array}
\end{displaymath}
For $m=3$, considering the bijection $F_3: \A^3 \rightarrow \R^4$ that maps  $\{\xi_2, \xi_{01}, \xi_{10}, \xi_{000}\}$ to $\{e_i\}_{i=1}^4$,  we get 
\begin{displaymath}
J_3 = 
\left(
\begin{array}{rrrr}
0 & 0 & 1 & 0 \\
2 & 1 & -1 & 0 
\end{array}
\right).
\end{displaymath}
For $m=4$, considering the bijection $F_4: \A^4 \rightarrow \R^8$ that maps
$\{\xi_{3}, \xi_{02}, \xi_{20}, \xi_{11}, \xi_{001}, \xi_{010}, \xi_{100}, \xi_{0000}\}$ to $\{e_i\}_{i=1}^8$, we get
\begin{displaymath}
J_4 = 
\left(
\begin{array}{rrrrrrrr}
0 & 0 & 0 & 0 & 0 & 0 & 1 & 0\\
0 & 0 & 0 & 1 & 0 & 0 & 0 & 0\\
0 & 0 & 2 & 0 & 0 & 1 & -1 & 0\\
1 & 0 & 0 & 0 & 0 & 0 & 0 & 0\\
0 & 6 & -6 & 0 & -1 & 2 & -1 & 0
\end{array}
\right).
\end{displaymath}

\smallskip

\noindent\textbf{Step 7.} Construction of the projections of the elements $\{\ell_{j}\}_{j=1}^n$  onto the orthogonal complement to the right ideal $\J_{a,b}$ (the main part of the series).

For every $i=1,\ldots,n$ we construct the projection of $\ell_i$ onto $\J_{a,b}^\perp$, which we denote by $\widetilde{\ell}_i$. 

Let $\ord(\ell_i)=m$. Then, actually, we  construct the projection of $\ell_i\in \A^{m}$ onto $\J_{a,b}^\perp \cap \A^{m}$. Let us consider the bijection $F_m:\A^{m}\to \R^{q}$ used at Step~6.  Since $\widetilde{\ell}_i\in \J_{a,b}^\perp$, we get  $J_{m} F_{m}(\widetilde{\ell}_i)=0$. 

Let us denote $y_i=\ell_i-\widetilde{\ell}_i\in \J_{a,b}\cap \A^{m}$. Then $F_{m}(y_i)$ equals a linear combination of the (transposed) rows of the matrix $J_{m}$, i.e., $F_{m}(y_i)=J_{m}^T x_i$ for some column vector $x_i$. On the other hand, $J_{m} F_m(y_i)=J_{m} F_m(\ell_i-\widetilde{\ell}_i)=J_{m} F_m(\ell_i)$.

Therefore, we get the following system of linear equations for $x_i$
\begin{equation}\label{eq_projection}
J_{m} J_{m}^T x_i = J_{m} F_m(\ell_i).
\end{equation}
The matrix  $J_{m}$ is of full rank, hence, the system \eqref{eq_projection} has a unique solution. Having found it, we then calculate $F_m(y_i)=J_{m}^T x_i$, i.e., express $y_i$ as a linear combination of basis elements of $\J_{a,b}\cap\A^m$ (collected in $D_j^m$). Then we find $\widetilde\ell_i=\ell_i-y_i$. 

For the system \eqref{eq_ex1} we have $\widetilde{\ell}_1 = \ell_1 = \xi_0$ since $\J_{a,b} \cap \A^{1} = \{0\}$. Further, for $\ell_2$,  using the bijection $F_3$ and the matrix $J_3$ from the previous step, we find $x_2$ as a solution of the the equation \eqref{eq_projection}, which takes the form
\begin{displaymath} 
\left(
\begin{array}{rr}
1 & -1 \\
-1 & 6 
\end{array}
\right) x_2 = 
\left(
\begin{array}{r}
0 \\
2
\end{array}
\right).
\end{displaymath}
Thus, $x_2 = \frac25 (1,1)^T$ and therefore $y_2$ equals $\frac25$ times the sum of the elements from $D_1^3$ and $D_2^3$, i.e., $y_2=\frac25(\xi_{10}+2\xi_2+\xi_{01}-\xi_{10})=\frac25 (\xi_{01}+2\xi_2)$. Hence, 
\begin{displaymath}
\widetilde{\ell}_2 = \ell_2 - y_2 = \xi_2-{\textstyle\frac25} (\xi_{01}+2\xi_2)={\textstyle\frac15} \xi_2 - {\textstyle\frac25}\xi_{01}.
\end{displaymath}
Applying the same procedure to $\ell_3$ we get
\begin{displaymath}
\widetilde{\ell}_3 = {\textstyle\frac{3}{19}}\xi_{02} + {\textstyle\frac{23}{285}}\xi_{20} +{\textstyle\frac{8}{57}}\xi_{001} - {\textstyle\frac{46}{285}}\xi_{010}.
\end{displaymath}

\smallskip

\noindent\textbf{Step 8.} Construction of the non-autonomous approximating system.

We follow the procedure described in Subsection~\ref{sec24}. First we implement recursively the shuffle product for basis elements:
\begin{displaymath}
\xi_{m_1\ldots m_k} \shf \xi_{n_1\ldots n_r}.
\end{displaymath}
The result is a homogeneous sum of basis elements and we represent it as a vector in the corresponding space $\R^q$, $q=2^{\sum m_i + \sum n_j + k + r -1}$. Then we extend by linearity the shuffle operation on arbitrary homogeneous elements of the algebra $\A$.

Further for every $i=1,\ldots,n$ we construct the function $\widehat{b}_i(t,x)$ based on the representation of $\widetilde{\ell}_i\in \A_{m}$:
\begin{displaymath}
\widetilde{\ell}_i = \sum_{j=0}^{m-2} y_j  \xi_j + \alpha_i \xi_{m-1}.
\end{displaymath}
We represent each $y_j$ as a shuffle polynomial with respect to $\widetilde{\ell}_1, \ldots, \widetilde{\ell}_{i-1}$. To find all shuffle products of $\{\widetilde{\ell}_{j}\}$ of the same order $m-j-1$, we use the method similar to one described at the Step~1.

Let us apply this procedure to the system \eqref{eq_ex1}. Since  $\widetilde{\ell}_1 = \xi_0$, we have $\widehat{b}_1(t,x)=-1$. Further,
\begin{displaymath}
\widetilde{\ell}_2 = {\textstyle\frac15} \xi_2- {\textstyle\frac25}\xi_{0}\xi_{1}  
= {\textstyle\frac15} \xi_2- {\textstyle\frac25}\widetilde{\ell}_1\xi_{1},
\end{displaymath}
which gives $\widehat{b}_2(t,x)= - \frac15 t^2+\frac25 t x_1$.
Finally,
\begin{displaymath}
\widetilde{\ell}_3 = {\textstyle\frac{3}{19}}\xi_{0}\xi_{2} + {\textstyle\frac{8}{57}}\underbrace{\xi_{00}}_{=\frac12\widetilde\ell_1^{\shf 2}}\xi_{1}+
{\textstyle\frac{23}{285}}(\underbrace{\xi_{2} - 2\xi_{01}}_{=5\widetilde\ell_2})\xi_{0}
\end{displaymath}
\begin{displaymath}
= {\textstyle\frac{3}{19}}\widetilde{\ell}_1\xi_{2} + {\textstyle\frac{4}{57}}\widetilde{\ell}_1^{\shf 2}\xi_{1}+ {\textstyle\frac{23}{57}}\widetilde{\ell}_2\xi_{0}, 
\end{displaymath}
which gives $\widehat{b}_3(t,x)=-\frac{3}{19}t^2 x_1-\frac{4}{57}t x_1^2 -\frac{23}{57}x_2$, and we obtain the following non-stationary approximation:
\begin{equation}\label{eq_ex1_appr}
\left\{
\begin{array}{l}
\dot{x}_1 = -u, \\
\dot{x}_2 = (- \frac15t^2+\frac25 t x_1) u, \\
\dot{x}_3 = (-\frac{3}{19}t^2 x_1 - \frac{4}{57} t x_1^2 - \frac{23}{57} x_2)u.
\end{array}
\right.
\end{equation}

\smallskip

\noindent\textbf{Step 9.} Construction of  the autonomous approximating system.

We follow the procedure described in Subsection~\ref{sec25}. To this end, we implement the linear maps $\varphi$ and $\psi$ in the algebra $\A$
and use the shuffle product implemented at the previous step.

For every $i=1,\ldots,n$ we find $\varphi(\widetilde{\ell}_i)$, $\psi(\widetilde{\ell}_i)$
and represent them as shuffle polynomials with respect to $\widetilde{\ell}_1, \ldots, \widetilde{\ell}_{i-1}$ if possible. It may happen that the corresponding system for $\varphi(\widetilde{\ell}_i)$ is unsolvable, which means that there is no autonomous approximating system. 
Again, to find all shuffle products of $\{\widetilde{\ell}_{j}\}$ of a given order, we use the method similar to one described at Step~1.

Applying this procedure to the system \eqref{eq_ex1}, we see that the autonomous approximating system does not exist since
$\varphi(\widetilde{\ell}_2) = \frac25(\xi_{1} - \xi_{00})$ cannot be represented as a shuffle polynomial with respect to $\widetilde{\ell}_1=\xi_0$.

Let us change the system \eqref{eq_ex1} a little bit by adding the term $-2tx_1$ to the second equation: 
\begin{equation}\label{eq_ex1_changed}
\left\{
\begin{array}{l}
\dot{x}_1 = - u \cos x_1,\\
\dot{x}_2 = -\sin^2 x_1 -2tx_1 + t^2 u, \\
\dot{x}_3 = 2x_1^2 \sin t - x_2 u.
\end{array}
\right.
\end{equation}
For this system, $v(\xi_2) = 0$ and 
\begin{displaymath}
\ell_{1} = \xi_0,\quad \ell_{2} = [\xi_0, \xi_1],\quad \ell_{3} = [\xi_0, \xi_2].
\end{displaymath}
Applying the described steps to the system \eqref{eq_ex1_changed}, we construct the following autonomous approximation:
\begin{displaymath}
\left\{
\begin{array}{l}
\dot{x}_1 = - u, \\
\dot{x}_2 = -\frac12 x_1^2, \\
\dot{x}_3 = \frac{1}{27}x_1^3 - \frac{10}{9} x_2 + \frac{4}{9} x_2 u.
\end{array}
\right.
\end{displaymath}

\section{Implementation of the algorithm and perspectives} 

The described above functionality and the pipeline of the algorithm were implemented as a package and a web application by using Python programming language. Here are the links: 

the web application: 
 \href{https://happycontrol.pythonanywhere.com/}{https://happycontrol.pythonanywhere.com/},

the package: 
 \href{https://github.com/ViktorRusakov/happy-control}{https://github.com/ViktorRusakov/happy-control}.

\begin{figure}[h]
	\centering
	\framebox{\includegraphics[width=0.75\textwidth]{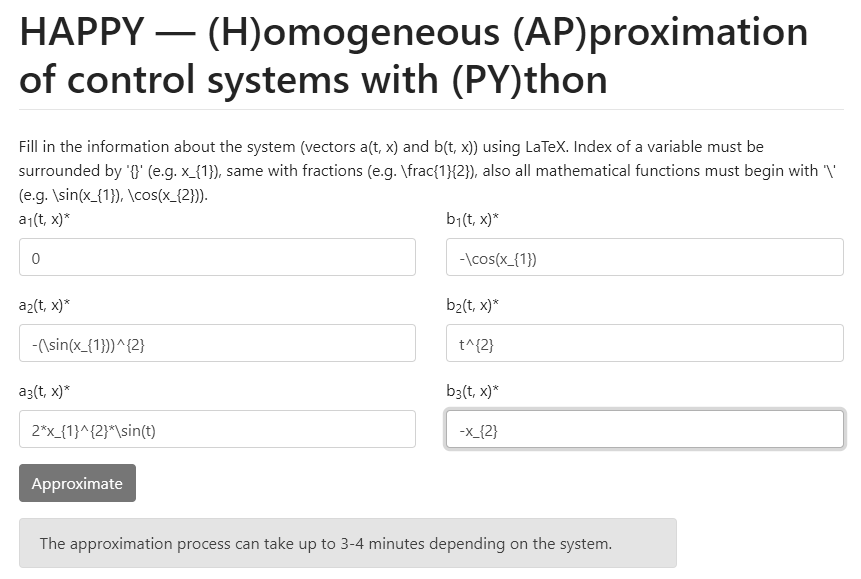}}
	\caption{Screenshot of the web application: setting a system for approximation}
	\label{fig1}
\end{figure}

The web application was built with the Django framework and has the following functionality. First, one has to choose system's dimension and type of the problem to be solved. Then one enters coefficients of the system by using LaTeX text formatting (and possibly the point in which neighborhood  the system is considered). The output is a pdf file with detailed calculations similar to the ones from the example given in the previous section. In Fig.~\ref{fig1} we give a screenshot of the application running.
The functionality on which the application is built, was implemented as an open source package, so it is possible to use it both as an end-to-end solution or as particular functions.

The core packages used for calculations are SymPy and NumPy:

\begin{itemize}
\item SymPy is used for symbolic calculations when constructing series representation of control system (Step~3 requires differentiation operation) and when printing out nice formatted output pdf files;

\item NumPy is used for numerical calculations: finding linearly independent Lie elements, their projections, etc. Basically all matrix operations were done with NumPy.
\end{itemize}

The most computationally demanding are Steps 7, 8 and 9 since they require solving high dimensional systems of linear equations. Since values of $\dim \A^m$ grow as $2^{m-1}$, one can conclude that the ``depth'' of a system, i.e., the order of the last linearly independent Lie element $\ord( \ell_n)$, determines the complexity and time of computations rather than dimension $n$ of a system.

The present-day implementation of the algorithm possesses the following restrictions on the depth of systems: the graded subspaces and the basis of the Lie algebra were constructed for orders up to $15$; the web application can handle in real time systems with ${\ord(\ell_n) \le 9}$; by using the functionality of the package one can construct approximations for systems with ${\ord(\ell_n) \le 13}$ on an ``average'' local computer in a ``reasonable'' time.

As possible developments of the proposed algorithm and its implementation we consider the following:

\begin{enumerate}

\item[1)] optimization of the present implementation to be able to construct approximations for ``deeper'' systems;  

\item[2)] implementation of the algorithm for the case of multi-dimensional control (for that, optimization is even more crucial since dimensions of graded subspaces for a system with $k$ controls grows as $(k+1)^{m-1}$);

\item[3)] applying the algorithm as an effective tool for analysis of certain specific classes of systems, in particular, Goursat systems, see, e.g., \cite{Mormul:06}.

\end{enumerate}

\section{Physical example: Goursat systems and their homogeneous approximations} 

Let us consider a driftless control system \eqref{eq:G} defined by two vector fields and introduce the following notation for the ``derived flag'' 
\begin{displaymath}
	\begin{array}{l}
	D^1=\Lin\{X_1,X_2\},\\ 
	D^{k+1}=D^{k}+[D^{k},D^{k}], \ k\ge1,
	\end{array}
\end{displaymath}
where $[D^{k},D^{k}]=\Lin\{[Y,Z]:Y,Z\in D^k\}$ and $[\cdot,\cdot]$ stands for the Lie bracket of vector fields. By definition, \eqref{eq:G} is called a Goursat system if  $\dim(D^k)=k+1$ for $k=1,\ldots,n-1$. A convenient coordinate description was found in \cite{Kumpera:80}: it turned out that, in appropriate local coordinates, vector fields $X_1(x)$ and $X_2(x)$ are polynomial and can include constant parameters. 

These systems have an important application discovered by F.~Jean \cite{Jean:96}: it relates to a kinematic model of a car pulling $n$ trailers. Roughly speaking, variables here are orientation angles of the car and the trailers. If the car is not perpendicular to the first trailer and none of the consecutive trailers are perpendicular to each other, then the position is non-singular. In this case, constant parameters vanish, and all non-singular positions are described by a unique system (of a chained form). Otherwise, different singularities arise; a thorough research was made in \cite{Jean:96}. 

The study of possible homogeneous approximations of Goursat systems, which describe a possible local behavior of the mentioned kinematic model, is very difficult in general; the main question here is if constant parameters mentioned above are included in a homogeneous approximation. A deep investigation was carried out by P.~Mormul \cite{Mormul:05,Mormul:06}, who performed a lot of sophisticated calculations applying very special properties of Goursat systems. Actually, Goursat systems were our main motivation in developing an implementation of the algorithm presented above. 

The described above restrictions on the depth of systems for the current implementation of the algorithm applies to Goursat systems as well. One of the ``deepest'' system for which we were able to calculate a homogeneous approximation is GGSGSG (written in a standard coding for Goursat systems, we refer to \cite{Jean:96}, \cite{Mormul:06} for more details): this is a system of dimension 8 and Lie brackets up to length 12 are involved. However, we believe that based on the special structure of Goursat systems, the implementation of the algorithm can be significantly optimized, and we will be able to find homogeneous approximations for Goursat systems of higher dimension.   

\subsection*{Acknowledgments}

The work was financially supported by Polish National Science Centre grant no. 2017/25/B /ST1/01892.
Pavel Barkhayev was supported by the Norwegian Research Council project ''COMAN'' No. 275113.

\bibliographystyle{plain}
\bibliography{BIRS_bib}

\begin{thebibliography}{10}

\bibitem{Agrachev_Gamkrelidze_Sarychev:89}
A.~A. Agrachev, R.~V. Gamkrelidze, and A.~V. Sarychev.
\newblock Local invariants of smooth control systems.
\newblock {\em Acta Appl. Math.}, 14:191--237, 1989.
\newblock DOI: 10.1007/BF01307214.

\bibitem{Bellaiche:96}
A.~Bella\"{\i}che.
\newblock The tangent space in sub-{R}iemannian geometry.
\newblock In A.~Bella\"{\i}che and J.~J. Risler, editors, {\em Sub-{R}iemannian
  {G}eometry}, volume 144 of {\em Progress in Mathematics}, pages 1--78.
  Birkh\"{a}user Basel, 1996.
\newblock DOI: 10.1007/978-3-0348-9210-0\_1.

\bibitem{Bianchini_Stefani:90}
R.~M. Bianchini and G.~Stefani.
\newblock Graded approximation and controllability along a trajectory.
\newblock {\em SIAM J. Control Optimiz.}, 28:903--924, 1990.
\newblock DOI: 10.1137/0328050.

\bibitem{Bressan:85}
A.~Bressan.
\newblock Local asymptotic approximation of nonlinear control systems.
\newblock {\em Internat. J. Control}, 41:1331--1336, 1985.

\bibitem{Chen}
K.~T. Chen.
\newblock Integration of paths -- a faithful representation of parths by
  noncommutative formal power series.
\newblock {\em Trans. Amer. Math. Soc}, 89:395--407, 1958.

\bibitem{Crouch:84}
P.~E. Crouch.
\newblock Solvable approximations to control systems.
\newblock {\em SIAM J. Control Optimiz.}, 22:40--54, 1984.
\newblock DOI: 10.1137/0322004.

\bibitem{Fliess:78}
M.~Fliess.
\newblock D\'{e}veloppements fonctionnels en ind\'{e}termin\'{e}es non
  commutatives des solutions d'\'{e}quations diff\'{e}rentielles non
  lin\'{e}aires forc\'{e}es.
\newblock {\em C. R. Acad. Sci. Paris Ser. A-B}, 287:1133--1135, 1978.

\bibitem{Fliess:80}
M.~Fliess.
\newblock Realizations of nonlinear systems and abstract transitive {L}ie
  algebras.
\newblock {\em Bull. Amer. Math. Soc.}, 2:444--446, 1980.
\newblock DOI: 10.1090/S0273-0979-1980-14760-6.

\bibitem{Fliess:81}
M.~Fliess.
\newblock Fonctionnelles causales non lin\'{e}aires et ind\'{e}termin\'{e}es
  non commutatives.
\newblock {\em Bull. Soc. Math. France}, 109:3--40, 1981.

\bibitem{Hermes:86}
H.~Hermes.
\newblock Nilpotent approximations of control systems and distributions.
\newblock {\em SIAM J. Control Optimiz.}, 24:731--736, 1986.
\newblock DOI: 10.1137/0324045.

\bibitem{Hermes:91}
H.~Hermes.
\newblock Nilpotent and high-order approximations of vector field systems.
\newblock {\em SIAM Rev.}, 33:238--264, 1991.
\newblock DOI: 10.1137/1033050.

\bibitem{Ignatovich:09}
S.~Yu. Ignatovich.
\newblock Realizable growth vectors of affine control systems.
\newblock {\em J. Dyn. Control Syst.}, 15:557--585, 2009.
\newblock 10.1007/s10883-009-9075-y.

\bibitem{Ignatovich:11}
S.~Yu. Ignatovich.
\newblock Normalization of homogeneous approximations of symmetric affine
  control systems with two controls.
\newblock {\em J. Dyn. Control Syst.}, 17:1--48, 2011.
\newblock 10.1007/s10883-011-9109-0.

\bibitem{Isidori}
A.~Isidori.
\newblock {\em Nonlinear Control Systems}.
\newblock Springer-Verlag London, 3rd edition, 1995.

\bibitem{Jean:96}
F.~Jean.
\newblock The car with $n$ trailers: {C}haracterization of the singular
  configurations.
\newblock {\em ESAIM Control Optim. Calc. Var.}, 1:241--266, 1996.
\newblock 10.1051/cocv:1996108.

\bibitem{Kawski:97}
M.~Kawski.
\newblock Nonlinear control and combinatorics of words.
\newblock In {\em Geometry of Feedback and Optimal Control}, pages 305--346.
  Dekker, 1997.

\bibitem{Kawski:01}
M.~Kawski.
\newblock The combinatorics of nonlinear controllability and noncommuting
  flows.
\newblock volume~8 of {\em Abdus Salam ICTP Lecture Notes Series}, pages
  223--312. 2002.

\bibitem{Kawski_Sussmann:97}
M.~Kawski and H.~J. Sussmann.
\newblock Noncommutative power series and formal {L}ie-algebraic techniques in
  nonlinear control theory.
\newblock In {\em Operators, Systems and Linear Algebra. European Consortium
  for Mathematics in Industry}, pages 111--128. Teubner, 1997.
\newblock DOI: 10.1007/978-3-663-09823-2\_10.

\bibitem{Korobov_Sklyar:87}
V.~I. Korobov and G.~M. Sklyar.
\newblock Time-optimality and the power moment problem ({R}ussian).
\newblock {\em Mat. Sb. (N.S.)}, 134(176)(2):186--206, 1987.
\newblock Translation: Math. USSR-Sb., vol. 62, no. 1, pp. 185-206, 1989, DOI:
  10.1070/SM1989v062n01ABEH003235.

\bibitem{Korobov_Sklyar:91}
V.~I. Korobov and G.~M. Sklyar.
\newblock The {M}arkov moment min-problem and time optimality ({R}ussian).
\newblock {\em Sibirsk. Mat. Zh.}, 32(1):60--71, 1991.
\newblock Translation: Siberian Math. J., vol. 32, no. 1, pp. 46-55, 1991, DOI:
  10.1007/BF00970159.

\bibitem{Kumpera:80}
A.~Kumpera and C.~Ruiz.
\newblock Sur l'\'{e}quilalence locale des syst\`{e}mes de {P}faff en drapeau.
\newblock In {\em Monge-Amp\`{e}re Equations and Related Topics, Florence
  1980}, Inst. Alta Math. F. Severi, Roma, pages 201--248. 1982.

\bibitem{Mormul:05}
P.~Mormul.
\newblock Real moduli in local classification of {G}oursat flags.
\newblock {\em Hokkaido Math J.}, 34:1--35, 2005.
\newblock DOI: 10.14492/hokmj/1285766199.

\bibitem{Mormul:06}
P.~Mormul.
\newblock Do moduli of {G}oursat distributions appear on the level of nilpotent
  approximations?
\newblock In J.P. Brasselet and M.A.S. Ruas, editors, {\em {R}eal and {C}omplex
  {S}ingularities}, volume~39 of {\em {T}rends in {M}athematics}, pages
  229--246. Birkh\"{a}user Basel, 2006.
\newblock DOI: 10.1007/978-3-7643-7776-2\_17.

\bibitem{Reutenauer}
C.~Reutenauer.
\newblock {\em Free {L}ie Algebras}.
\newblock Clarendon Press, Oxford, 1993.

\bibitem{Sklyar_Ignatovich:96ZAMM}
G.~M. Sklyar and S.~Yu. Ignatovich.
\newblock Local equivalence of time-optimal control problems and {M}arkov
  moment min-problem.
\newblock {\em Z. Angew. Math. Mech.}, 76(S.3):561--562, 1996.

\bibitem{Sklyar_Ignatovich:00}
G.~M. Sklyar and S.~Yu. Ignatovich.
\newblock Moment approach to nonlinear time optimality.
\newblock {\em SIAM J. Control Optimiz.}, 38(6):1707--1728, 2000.
\newblock DOI: 10.1137/S0363012997329767.

\bibitem{Sklyar_Ignatovich:02}
G.~M. Sklyar and S.~Yu. Ignatovich.
\newblock Representations of control systems in the {F}liess algebra and in the
  algebra of nonlinear power moments.
\newblock {\em Systems Control Lett.}, 47(3):227--235, 2002.
\newblock DOI: 10.1016/S0167-6911(02)00201-3.

\bibitem{Sklyar_Ignatovich:03}
G.~M. Sklyar and S.~Yu. Ignatovich.
\newblock Approximation of time-optimal control problems via nonlinear power
  moment min-problems.
\newblock {\em SIAM J. Control Optimiz.}, 42:1325--1346, 2003.
\newblock DOI: 10.1137/S0363012901398253.

\bibitem{Sklyar_Ignatovich:07}
G.~M. Sklyar and S.~Yu. Ignatovich.
\newblock Description of all privileged coordinates in the homogeneous
  approximation problem for nonlinear control systems.
\newblock {\em C. R. Math. Acad. Sci. Paris}, 344:109--114, 2007.
\newblock DOI: 10.1016/j.crma.2006.11.016.

\bibitem{Sklyar_Ignatovich:08}
G.~M. Sklyar and S.~Yu. Ignatovich.
\newblock Fliess series, a generalization of the {R}ee's theorem, and an
  algebraic approach to a homogeneous approximation problem.
\newblock {\em Int. J. Control}, 81:369--378, 2008.
\newblock DOI: 10.1080/00207170701561427.

\bibitem{Sklyar_Ignatovich:14}
G.~M. Sklyar and S.~Yu. Ignatovich.
\newblock Free algebras and noncommutative power series in the analysis of
  nonlinear control systems: an application to approximation problems.
\newblock {\em Dissertationes Math. (Rozprawy Mat.)}, 504:1--88, 2014.
\newblock DOI: 10.4064/dm504-0-1.

\bibitem{Sklyar_Ignatovich:20_conf}
G.~M. Sklyar and S.~Yu. Ignatovich.
\newblock Constructing of a homogeneous approximation.
\newblock In A.~Bartoszewicz, J.~Kabzi\'{n}ski, and J.~Kacprzyk, editors, {\em
  Advanced, Contemporary Control}, volume 1196 of {\em Advances in Intelligent
  Systems and Computing}, pages 611--624. Springer, Cham., 2020.
\newblock DOI: 10.1007/978-3-030-50936-1\_52.

\bibitem{Sklyar_Ignatovich:21}
G.~M. Sklyar and S.~Yu. Ignatovich.
\newblock Subspaces of maximal singularity for homogeneous control systems.
\newblock {\em J. Dyn. Control Syst.}, 27:585--616, 2021.
\newblock DOI: 10.1007/s10883-020-09518-x.

\bibitem{Sklyar_Ignatovich_Barkhayev:05}
G.~M. Sklyar, S.~Yu. Ignatovich, and P.~Yu. Barkhayev.
\newblock Algebraic classification of nonlinear steering problems with
  constraints on control.
\newblock In {\em Advances in mathematics research. {V}ol. 6}, pages 37--96.
  Nova Sci. Publ., New York, 2005.

\bibitem{Stefani:85}
G.~Stefani.
\newblock Polynomial approximations to control systems and local
  controllability.
\newblock In {\em 1985 24th IEEE Conference on Decision and Control}, pages
  33--38, 1985.
\newblock DOI: 10.1109/CDC.1985.268467.

\bibitem{Sussmann:92}
H.~J. Sussmann.
\newblock New differential geometric methods in nonholonomic path finding.
\newblock In {\em Systems, Models and Feedback: Theory and Applications},
  volume~12 of {\em Progress in Systems and Control Theory}, pages 365--384.
  Birkh\"{a}user, Boston, MA, 1992.
\newblock DOI: 10.1007/978-1-4757-2204-8\_24.

\bibitem{Zuyev:19}
A.~Zuyev and V.~Grushkovskaya.
\newblock On stabilization of nonlinear systems with drift by time-varying
  feedback laws.
\newblock In {\em 2019 12th International Workshop on Robot Motion and Control
  (RoMoCo), Poznan, Poland}, pages 9--14, 2019.
\newblock doi: 10.1109/RoMoCo.2019.8787353.

\end{thebibliography}

\end{document}